# ENERJİ ŞEBEKELERİNİN KARIŞIK TAM SAYILI DOĞRUSAL OLMAYAN PROGRAMLAMAYLA OPTİMUM TASARIMI: PROBLEMLER, OLASI ÇÖZÜMLER VE MODEL FORMÜLASYONU

*Handan AKÜLKER* [*]

*Erdal AYDIN*[**,***]

**Özet:** Enerji sektörü gelişen teknolojiyle ve enerjiye duyulan talebin artmasıyla dünya genelinde öncelikli alan haline gelmiştir. Günümüz enerji kaynağının çoğunu oluşturan fosil yakıtlardan çıkan sera gazları küresel ısınmaya sebep olmaktadır. Ek olarak, Türkiye'nin temel enerji kaynakları bakımından dışa bağımlı oluşu ve üretilen enerjinin optimum dağıtılmaması, enerji şebekelerinin optimum tasarımı ve şebekelerin uzun vadeli üretim planlaması ihtiyacını arttırmıştır. Mikro şebeke; dağıtık güç jeneratörleri, depolama birimleri ve yüklerden oluşan bir enerji şebekesidir. Bu çalışmada, karma tam sayılı doğrusal olmayan programlama (MINLP) ile fosil ve yenilenebilir enerji kaynakları içeren mikro şebekelerin optimum tasarımı ve üretim planlaması incelenmiştir. Fazladan elektrik üretimi, $CO_2$ emisyonu, su kıtlığı gibi problemlere yönelik olası çözümler ele alınmıştır. Son olarak, optimizasyon modeli oluşturulurken dikkate alınması gereken kısıtlar ve her bir ekipmana ait model denklemleri incelenmiştir.

**Anahtar Kelimeler:** Karışık (karma) tam sayılı doğrusal olmayan programlama (MINLP), Enerji şebekesi optimizasyonu, Mikro şebeke, Fazladan üretilen elektrik, Karbondioksit emisyonu, Su kıtlığı

## Mixed Integer Nonlinear Programming for Optimal Design of Energy Grids: Optimization Problems, Possible Solutions, and Model Formulation

**Abstract:** The energy sector has become priority around the world with developing technology and increasing power and energy demand. That all sources for energy production are not renewable increases greenhouse gas emissions and causes global warming. Turkey's dependence on foreign sources in terms of vital energy resources and the suboptimality of the energy distribution has increased the need for optimum design of energy grids and long-term production planning of established grids. Microgrid is an energy grid consisting of distributed generators, storage units and loads. In this study, optimum design and generation

---

[*] Boğaziçi Üniversitesi, Mühendislik Fakültesi, Kimya Mühendisliği Bölümü, 34342, Bebek, İstanbul
[**] Koç Üniversitesi, Mühendislik Fakültesi, Kimya ve Biyoloji Müh. Bölümü, 34450, Sarıyer, İstanbul
[***] Koç Üniversitesi-Tüpraş Enerji Merkezi (KUTEM), 33450, Sarıyer, İstanbul

İletişim Yazarı: Erdal AYDIN (eaydin@ku.edu.tr)

planning of microgrids using mixed integer nonlinear programming (MINLP) is investigated. The optimal design under conditions such as excess electricity, CO$_2$ emission, and water shortage are examined. Finally, the constraints of the optimization models and the empirical equations for each equipment are investigated.

## 1. GİRİŞ

Enerji fiyatlarındaki hızlı artış, iklim değişikliği ve teknolojik gelişmeler nedeniyle, enerji endüstrisi dünya çapında öncelikli bir alan olmuştur. Hızla artan küresel ısınma sorunu, dünyamızın geleceği için dikkate alınması gereken bir diğer önemli konudur. Fosil yakıtlar enerji üretimi için kullanıldığında, küresel ısınmaya neden olan sera gazları gibi yan ürünler oluşur. Bu nedenle; son yıllarda enerji üretiminde sera gazı salınımı yapmayan yenilenebilir enerji kaynakları üzerinde durulmaktadır (Abo-Elyousr ve Elnozahy, 2018; Shi ve diğ., 2019; Talebi ve diğ., 2016). Mikro şebeke; dağıtık güç jeneratörleri (DG), depolama birimleri ve yükleri içeren bir güç ağıdır. Bir mikro şebeke, belirli bir yerleşim yerinin enerji talebini karşılayabilir veya ana şebekenin talep karşılama yükünü azaltmak için ana şebekeye (ulusal şebeke) hizmet edebilir. Mikro şebekenin en gerekli parçası olan DG'ler; mikro türbinler, yakıt hücreleri, dizel motorlar, gaz türbinleri gibi kombine ısı ve güç (CHP) üretim sistemleri ile fotovoltaik hücreler ve rüzgar türbinleri gibi yenilenebilir enerji jeneratörlerinden meydana gelir (Aghaei ve Alizadeh, 2013; Tenfen ve Finardi, 2015).

Türkiye'de, enerji üretiminde kullanılan doğal gaz ve petrol gibi büyük ölçüde yenilenemeyen enerji kaynakları yönünden ciddi bir ithalat ihtiyacı bulunmaktadır. Bununla birlikte, elektrik iletim ve dağıtım şebekeleri arasında teknik ve teknik olmayan sebeplere bağlı kayda değer enerji verimsizlikleri mevcuttur. Yapılan çalışmalarda, bir şebekenin uzun dönemli üretim planlaması yapmasının teknik sorunları çözmede etkili olduğu gösterilmiştir. Enerji üretimi planlanırken, talep tarafı yönetimi ve dağıtık üretim jeneratörlerinin birlikte düşünülmesi gerektiği sonucuna ulaşılmıştır (Babacan ve Unvan, 2020; Onat, 2010). Bu nedenle, bir elektrik üretim şebekesinin optimum tasarımı ve uzun vadeli optimum üretiminin planlanması hem Türkiye'nin enerji açığını hem de enerji kaybını büyük ölçüde azaltabilir.

Enerji şebeke sistemlerinin MINLP kullanılarak optimum şekilde tasarlanması ve işletilmesi konusunda yapılan çalışmalar son yıllarda artmaya başlamıştır. Honarmand ve diğ. (2021) çalışmalarında, yenilenebilir enerji kaynaklarına ve üç adet elektrik, soğutma ve termal depolama sistemine sahip bir enerji merkezinin işletme problemini çözmek için sağlam optimizasyon modeli önermişlerdir. Çalışma merkezi, İran'ın Hamedan kentindeki bir hastane binasında gerçek bir vaka çalışmasıdır. MINLP problemini DICOPT kullanarak GAMS'te çözmüşlerdir. Bir başka çalışmada, stokastik talepler ve yenilenebilir kaynaklar varlığında adalı mikro şebekelerin optimum çalışması için stokastik MINLP modeli kurulmuştur. Önerilen formülasyonda, MINLP problemi doğrusallaştırılmış ve problem AMPL'de CPLEX kullanılarak çözülmüştür (Vergara ve diğ., 2020). Bir diğer çalışmada, pratik kısıtlara sahip yerleşim binalarının enerji talebini karşılamak için en iyi teknoloji kombinasyonunun belirlenmesi amaçlanmıştır. Bu bölge planlaması, yenilenebilir ve temiz enerjinin uygulanması için akıllı bir mikro şebeke kurmayı amaçlamaktadır. Gaz türbini, absorpsiyonlu soğutma grubu, elektrikli soğutma grubu, yoğuşmalı kazan, toprak kaynaklı ısı pompası, PV, elektrokimyasal depolama, ısı depolama, buz depolama klima sistemi vb. gibi bir dizi teknoloji, aday ekipmanlar olarak belirlenmiştir. Bu çok amaçlı optimizasyon problemini çözmek için bir MINLP modeli geliştirilmiştir. Model GAMS'te Lindo çözücüsü kullanılarak çözülmüştür (Zheng ve diğ., 2018). Bir başka optimum tasarım çalışması da Montoya ve diğ. (2020) tarafından gerçekleştirilmiştir. Çalışmada, dağıtık güç jeneratörlerinin optimum lokasyonu ve boyutları GAMS kullanılarak BONMIN çözücüsü ile hesaplanmıştır.

Elektrik enerjisi üretim planlaması yapılırken alternatif akım (AC) güç akışı kısıtının dahil edildiği çalışmalar da literatürde mevcuttur. Örneğin; Santos ve diğ. (2020) pil enerji depolama sistemi yerleştirme problemini ele almak için MINLP modeli önermişlerdir. Modelde, AC güç akışı ve voltaj büyüklüğü kısıtlarını dikkate almışlardır. Başka bir çalışmada, Shuai ve diğ. (2020) bir mikro şebekenin gerçek zamanlı işletim problemini ele almışlardır. AC güç akışı ve pillerin operasyonel karakter özelliklerini probleme dahil ederek stokastik MINLP modeli geliştirmişlerdir.

Son dönemde popüler hale gelen konut tipi mikro şebekelerin (r-mikro şebeke) termal konfor seviyesini ve ekonomiyi temel alan optimum planlaması MINLP kullanılarak yapılmıştır. Dal-budak algoritması kullanılan problem çözümü LINGO11 programında yapılmıştır (Liu ve diğ., 2018). Başka bir ısıtma-soğutma planlaması araştırması Deng ve diğ. (2017) tarafından yapılmıştır. Çalışmada, Tianjin'deki bir enerji istasyonunun elektrikli soğutma grubu sistemi, toprak kaynaklı ısı pompası sistemi, su termal enerji depolama sistemi ile entegre bir enerji istasyonunun günlük işletme maliyetini en aza indirmek için MINLP ile işletme temelli optimum zamanlama stratejisi önerilmektedir.

Enerji üretim teknolojilerinin gelişmesi ve çeşitlenmesiyle birlikte; geleceğin enerji şebekeleri, çok taşıyıcılı (kaynaklı) mikro şebekelerden oluşacaktır. Amir ve diğ. (2019) çalışmalarında, güvenilirlik kriterlerini uygulayan ana şebekeye bağlı çok taşıyıcılı mikro şebekeyi ele almıştır. Şebeke bileşenlerin optimum tipini, boyutunu ve optimum sevkiyatı bulmak için bir MINLP modeli önermişlerdir. Çok taşıyıcılı mikro şebeke örneklerinin en güncel olanlarından biri de elektrik ve doğal gazın entegre edildiği sistemlerdir. Bu sistemler, özellikle yüksek verimlilik ve ekonomik hususlar nedeniyle son yıllarda büyük ilgi görmektedir. Mansouri ve diğ. (2020) çalışmalarında, ilk aşamada parçacık sürüsü optimizasyonu (PSO) algoritması kullanarak bir enerji merkezinin tasarımını gerçekleştirmişlerdir. İkinci aşamada ise MINLP ile ilk aşamada optimum şekilde tasarlanan bu enerji merkezinin işletme problemini çözmüşlerdir. Entegre enerji şebekelerinde istatistiksel işletme modelleri oluşturulmasıyla ilgili başka bir çalışma Alipour ve diğ. (2018) tarafından gerçekleştirilmiştir. Bu çalışmada, talep-yanıt programı ile MINLP kullanılmıştır. Bir diğer çalışmada, hibrit bir mikro şebekenin 24 saatlik çalışma düzeninin oluşturulması için stokastik MINLP modeli geliştirilmiştir. Çalışmada yük tahmin hatasının etkisi ve enerji depolamanın kullanılabilirliği gibi etmenler incelenmiştir (Alvarado-Barrios ve diğ., 2020).

Su kıtlığı probleminin ortaya çıkmasıyla, enerji şebekesi optimizasyonu çalışmaları, su-enerji entegrasyonlu mikro şebeke çalışmaları alanına doğru kaymaya başlamıştır. Su ve enerji talepleri karşılanırken, bu iki önemli ihtiyacın ayrı ayrı düşünülmemesi gerektiği yapılan son çalışmalarda detaylı şekilde gösterilmiştir. Su dağıtım sistemleri, tıpkı su arıtma sistemlerinde olduğu gibi elektrik tüketimine ihtiyaç duyar. Örneğin; Li ve diğ. (2019) çalışmalarında, dağıtım düzeyinde bir su-enerji bağlantısı çerçevesinde yenilenebilir kaynakların yönetimi için, su şebekesinin elektrik şebekesine talep yanıt potansiyelini araştırmıştır. Geliştirdikleri MINLP ile talep tarafı yönetimi optimize edilmiştir. Öte yandan, Moazeni ve Khazaei (2021) çalışmalarında bir su-enerji şebekesine tersine pompa eklemenin etkisini incelemiştir. Enerji şebekesi kısmı, doğal gaz ve dizel bazlı jeneratörler, güneş fotovoltaik üniteleri, rüzgar jeneratörleri ve pil enerji depolama sistemlerinden oluşmaktadır. Geliştirilen MINLP modeli ile konulması düşünülen tersine pompaların sayısı, yeri ve üretilen enerjinin maliyeti optimize edilmiştir. Başka bir çalışmada, yine MINLP modeli ile su talebi ve enerji tüketimi birlikte optimize edilmiştir. Çalışmada adalı modda bir mikro su-enerji şebekesi incelenmiş ve bu şebekenin enerji verimliliği maksimize edilmek istenmiştir (Moazeni ve diğ., 2020). Aynı çalışma grubunun başka bir çalışmasında ise, yine adalı modda bir mikro su-enerji şebekesi incelenmiş; fakat bu çalışmada şebeke aynı zamanda akıllı bina termal sistemi ile entegre edilmiştir. 7 düğümlü bir su dağıtım sistemi, 4 akıllı bina ve hibrit bir mikro şebeke birlikte ele alınmış ve kurulan MINLP modeli ile dinamik ekonomik dağıtım hesaplanmıştır (Moazeni ve Khazaei, 2020).

Bu derleme çalışmasında, Karışık (Karma) Tam Sayılı Doğrusal Olmayan Programlama (MINLP) ile hem yenilenebilir hem de geleneksel jeneratörler içeren hibrit bir mikro şebekenin optimum tasarımı ve üretim planlaması incelenmiştir. Optimizasyonda karşılaşılan sorunlar, bu sorunların olası çözümleri ve model formülasyonları detaylı şekilde açıklanmıştır. Enerji sistemleri optimizasyonu çözüm kolaylığı sunması sebebiyle, çoğunlukla doğrusal ele alınır (Farrokhifar ve diğ., 2020; Feng ve diğ., 2019; Ren ve Gao, 2010). Doğrusal olmayan optimizasyon çalışmaları son yıllarda artmaya başlasa da, literatürdeki doğrusal çalışmalara kıyasla çok daha azdır. Problemin boyutunu ve karmaşıklığını azaltmak için doğrusal olmayan denklemler yerine doğrusal formların kullanılması, yakınsama sebebiyle karı ve maliyeti değiştirebilir. Optimum tasarımdaki ekipman sayısını ve çeşitliliğini de etkileyebilir ve sonuç olarak optimum olmayan bir sonuç ile karşılaşılabilir.

Mikro şebeke optimizasyonunda; şebeke bileşenlerinin ömürleri, sosyal veya düzenleyici kısıtlar, fazladan üretilen elektrik, çevresel ve atmosferik kısıtlar, $CO_2$ salınım limiti ve su kaynağı kıtlığı gibi birçok problemle karşılaşılmaktadır. Su arzının, küresel iklim sorunları sebebiyle enerji üretim sistemleri ile simbiyotik birleşim kurması son zamanlarda önemli bir kavram hâline gelmiştir. Bu çalışmada, fazladan elektrik üretimi, $CO_2$ emisyonu ve su temini problemleri ele alınmakta ve bu problemlerin çözümüne yönelik alternatifler sunulmaktadır. Ek olarak, genel hatlarıyla optimizasyon modelinin nasıl oluşturulacağı, MINLP çözüm algoritmaları ve çözücü tipleri hakkında bilgi verilmektedir.

## 2. MİKRO ŞEBEKE BİLEŞENLERİ

### 2.1. Kombine Isı ve Güç (CHP) Jeneratörleri

Kombine ısı ve güç (CHP) jeneratörleri, bir yakıt kaynağından aynı anda elektrik enerjisi ve termal enerji üretmek için kullanılırlar. Hem kullanıcının termal gereksinimlerini hem de yerleşimin elektrik talebinin tamamını veya bir kısmını karşılamak için, kullanıcının yerleşiminde veya yerleşimin yakınlarında kurulurlar (Aghaei ve Alizadeh, 2013). Çifte faydaları nedeniyle sera gazı salınımlarını azaltarak enerji verimliliğini artırırlar. CHP sistemleri, belirli parçalarla koordineli bir çerçevede kurulurlar. Bu parçalar temel olarak ana (primer) taşıyıcı birimi, jeneratör, ısı geri kazanım birimi ve elektrik ara bağlantısıdır. Ana taşıyıcı birimi, tüm sistemi çalıştırır ve sistemin tipi bu ana taşıyıcı birimine göre belirlenir. Başlıca CHP sistemleri, pistonlu içten yanmalı motorlar, yanmalı (gaz) türbinler, buhar türbinleri, mikro türbinler ve yakıt hücreleridir.

Her tip CHP sistemi aynı çalışma mekanizmasına sahiptir. Temelde kimyasal yakıtı elektrik enerjisine dönüştürürler; ancak bunu doğrudan yapmazlar. Kimyasal enerjiyi önce ısıya dönüştürürler, bu ısıyı elektrik üretmek için kullanılırlar. Yakıt hücreleri dışındaki tüm CHP sistemleri yakıtı yakarak ısı üretirler. Yakıt pilleri, kimyasal yakıtı elektrokimyasal bir yolla elektrik üretmek için kullanırlar. Kombine ısı ve güç jeneratörlerinin 5 farklı tipinin kullanım avantajları ve dezavantajları açısından karşılaştırılmıştır ve tablo şeklinde sunulmuştur:

**Tablo 1. Kombine ısı ve güç sistemlerinin avantajları ve dezavantajları** (EPA, 2017)

| Kombine Isı ve Güç Sistemleri | Avantajları | Dezavantajları |
|---|---|---|
| Pistonlu İçten Yanmalı Motorlar | *Kısmi yük işletme esnekliği ile yüksek güç verimliliği<br>*Cihazın hızlı açılması<br>*Diğer tiplere oranla düşük yatırım maliyeti<br>*Normal operatörler ile sahada revize edilebilirlik. | *Yüksek bakım maliyeti<br>*Düşük sıcaklıklarda sınırlanması<br>*Diğer tiplere oranla yüksek gaz salınımı<br>*Geri kazanılan ısıdan yararlanılmasa bile soğutma ihtiyacı. |
| Buhar türbinleri | *Yüksek toplam verimlilik<br>*Birden fazla yerin ısı gereksinimini karşılama yeteneği<br>*Uzun ömürlülük | *Cihazın yavaş açılması<br>*Çok düşük güç/ısı oranı<br>*Başka kazanlara ve buhar kaynaklarına duyulan ihtiyaç |
| Gaz türbinleri | *Yüksek dayanıklılık<br>*Düşük gaz salınımları<br>*Yüksek kalite ısı<br>*Soğutma ihtiyacı olmaması | *Yüksek basınçlı gazlara ve gaz kompresörlerine ihtiyaç duyulması<br>*Düşük yüklemede zayıf verimlilik<br>*Çıktının çevre sıcaklığı arttıkça düşmesi |
| Mikrotürbinler | *Az sayıda hareketli parça<br>*Kompakt boyut ve hafiflik<br>*Düşük gaz emisyonu<br>*Soğutma ihtiyacı olmaması | *Yüksek maliyet<br>*Diğerlerine oranla düşük mekanik verimlilik<br>*Düşük sıcaklıklarla sınırlı |
| Yakıt Hücreleri | *Düşük gaz salınımı<br>*Yüksek verimlilik<br>*Modüler tasarım | *Yüksek maliyet<br>*Yakıtın safsızlığına duyarlı<br>*Düşük güç yoğunluğu |

### 2.2. Fotovoltaik Hücreler

Fotovoltaik hücreler, doğrudan veya dolaylı olarak güneş radyasyonunu enerjiye dönüştürülebilen jeneratörlerdir. Güneş radyasyonu, beyaz ışık olarak hücreler tarafından alınır ve kızılötesinden ultraviyoleye kadar geniş bir dalga boyu spektrumuna yayılır. Bu radyasyonla, fotovoltaik hücre içinde doğrudan elektrik üretilebilmektedir. Fotovoltaik hücreler, tasarım ve kurulum açısından en kolay teknoloji olarak bilinmektedir; ancak yine de en pahalı yenilenebilir teknolojilerden biri olmayı sürdürmektedir.

Fotovoltaik hücreler tarafından üretilen maksimum güç, çoğunlukla bu hücrelerin imalatında kullanılan malzemelere bağlıdır. İlk nesil fotovoltaik hücre, kristal yapılı silikondan üretilmiştir. Kristal yapı, yıllar içinde geliştirilerek farklı sınıflara ayrılmıştır. Mono-kristalin ve poli-kristalin en bilinen silikon kristal yapı malzemeleridir (El Chaar ve diğ., 2011). Kristal silikon teknolojisini takip eden ikinci malzeme ince film teknolojisidir. Fotovoltaik hücrelerden en verimli şekilde elektrik üretebilmek için, koşullara uygun teknolojilerin kullanılması gerekmektedir. Hangi tip malzemeden üretilen fotovoltaik hücrelerin kullanılması gerektiği, hücrelerin kurulacakları iklim koşullarına göre yapılmalıdır. Yapılan bir çalışmada, kristal yapılı fotovoltaik hücreler güneşli günlerde bulutlu günlere oranla daha iyi performans göstermiştir. Öte yandan, aynı çalışmada ince film teknolojisiyle üretilen fotovoltaik hücrelerin bulutlu havalarda daha verimli olduğu görülmüştür (Premalatha ve Rahim, 2017).

## 2.3. Rüzgar Türbinleri

Rüzgar enerjisi yenilenebilir ve tükenmez bir enerji kaynağı olduğu için birçok ülke tarafından tercih edilmektedir. Rüzgar enerjisi, özellikle son 20 yılda dünya çapında önemli bir büyüme göstermiştir. Daha büyük boyutta rüzgar türbinlerinin kullanılması ve bu türbinlerin düzenli olarak yerleştirilmesi ile çiftlikler oluşturulması planlanmaktadır (Novaes Menezes ve diğ., 2018).

Daha büyük rüzgar türbinleri kurulurken bazı boyut sınırlamaları dikkate alınmalıdır. Güç çıktısı, türbin kanadı alanıyla doğru orantılıdır. Diğer bir deyişle, rüzgar türbini rotor kanadının çapı arttıkça, büyüyen rotor taramalı alan sayesinde enerji çıktısı artar. Ancak, artan alan nedeniyle artan hacim, çapın küpü kadar maliyeti yükseltmektedir.

Rüzgar türbinleri, rüzgar hızına göre belirlenen dört özel çalışma bölgesine sahiptir. İlk bölgede, elektrik enerjisi üretilmez çünkü rüzgar hızı 'cut-in' hızının altındadır. 'cut-in' rüzgar hızı, elektrik enerjisi üretiminin başladığı rüzgar hızı olarak tanımlanır. Rüzgar hızı cut-in hızının altındaysa rüzgar türbini elektrik tüketmeye başlayacaktır. Bu durumda türbinin kapatılması önerilir. İkinci bölge, cut-in ve nominal rüzgar hızı arasındaki rüzgar hızı bölgesi olarak tanımlanır. Rüzgar hızı kademeli olarak artarken enerji üretimi de artar. Rüzgar türbini bu bölgede kısmi yük rejimindedir. Üçüncü bölgede, rüzgar nominal hıza ulaşır. Rüzgar türbini, nominal güç üretir ve tam yük rejiminde çalışır. Rüzgar türbininin mekanik güvenliğini karşılamak için, güç üretimi rüzgar türbininin nominal gücü ile sınırlandırılmalıdır. Son bölge, 'cut-out' rüzgar hızından sonraki bölge olarak tanımlanır. 'cut-out' rüzgar hızı, rüzgar türbininin kapatılması gereken güvenlik sınır hızıdır (Amano, 2017).

## 2.4. Depolama Birimleri

Enerji depolama, literatürde nispeten yeni bir araştırma alanıdır ve birçok depolama tekniği henüz tam olarak geliştirilememiştir. Depolama teknolojilerinin avantajlarını ve dezavantajlarını anlamak, hangi uygulama için hangi depolama teknolojisinin en uygun olduğuna karar vermek gerekir. Enerji depolama birimleri dört temel kategoriye ayrılmaktadır: elektrokimyasal, mekanik, kimyasal ve termal. Bunlara ek olarak, enerji sistemi optimizasyonu için su depolama teknolojileri de dikkate alınmalıdır. Su depolama sistemleri enerji sistemiyle yakından ilgilidir; ancak doğrudan bir alt kategori olarak tanımlanmamıştır.

Jenerik isimleri 'pil' olan elektrokimyasal depolama teknikleri, enerji sistemleri tarafından üretilen fazla elektriği depolamak için kullanılır. Taşınabilir oldukları ve sistemdeki diğer ekipmanlara göre nispeten küçük oldukları için enerji sisteminin herhangi bir yerine yerleştirilebilirler. Şarj edilemez veya şarj edilebilir tipleri vardır. Lityum iyon, sodyum sülfür, kurşun asit ve redoks akışlı piller, en yaygın pil türleridir.

Mekanik enerji depolama, kinetik enerjinin fiziksel olarak elektrik enerjisine dönüştürülmesiyle gerçekleştirilir. Mekanik enerji depolama için iki önemli teknik vardır: basınçlı hava ve pompalı hidro depolama. Her iki teknik de üretim ve depolama yerlerinin konumuyla kısıtlanmıştır. Bu durum, çok uzun mesafeler için mekanik enerjinin depolanmanın zorlaşmasına sebep olmaktadır.

Kimyasal depolama, enerji sistemleri için çok önemli kimyasallar olan hidrojen ve metan depolamayı içeren önemli bir depolama teknolojisidir. Kimyasal enerji depolama mekanizması, kimyasal bağların sahip olduğu enerjiye dayanmaktadır. Endotermik kimyasal reaksiyonlar, yüksek enerjili kimyasal bağlar oluşturmak için enerjiye ihtiyaç duyarken, ekzotermik reaksiyonlar, daha düşük enerjili bağlar oluşturarak enerjiyi serbest bırakır. Bu durum, kimyasal bileşiklerin bağlarında elektrik ve ısıyı depolayan sistemlerin, bunları ilerde bir enerji kaynağı olarak tüketmesine olanak sağlar.

Düşük yoğunluklu hidrojen depolama, kaplara sıkıştırılarak veya sıvılaştırılarak gerçekleştirilir. Yüksek yoğunluklu depolama nanotüpler veya katı metal hidritler ile gerçekleştirilir. Yer altı tuz mağaraları, uzun vadede büyük miktarda hidrojen depolamak için kullanılabilir.

Metanı altyapı inşa ederek depolama ve nakliye için pratiktir. Doğal gaz şebekesinde de metan depolamak mümkündür (Gustavsson, 2016). Doğal gaz depolama, enerji sektörü için son derece önemli olmuştur. Bunun birinci nedeni doğal gazın yetersiz olması ve bu durumun enerji üreticileri arasındaki rekabeti canlandırmasıdır. İkinci neden, doğal gaz depolamanın diğer bileşenlerin yanı sıra esas olarak hidrokarbonlardan oluşan içeriği nedeniyle çevre dostu bir teknoloji olmasıdır. Doğal gaz genellikle yer altında büyük depolama rezervuarlarında depolanır. Doğal gazın depolanması için en çok bilinen diğer yöntem sıvılaştırmadır (LNG). Bu yöntemde doğalgaz atmosferik basınçta -160 ºC civarında soğutulur. Buharlaştığı zaman hava ile karıştırıldığında sadece birkaç farklı yoğunlukta yanar. Patlama riski olmadığından hem LNG hem de buharı çevreye zararsızdır (Eren ve Polat, 2020; Zhichao ve diğ., 2015).

Termal enerji depolama teknolojisi, bir depolama ortamını ısıtarak veya soğutarak termal enerjiyi depolar. Ortam içinde depolanan enerji daha sonra soğutma, ısıtma ve güç üretimi uygulamaları için tüketilebilir. Termal enerji depolaması için üç teknik vardır. Birincisi, depolamanın malzemedeki sıcaklık farkına dayandığı gelişmiş bir teknolojidir. Öte yandan, depolama ortamının özgül ısısına bağlı sınırlı bir depolama kapasitesine sahiptir. İkinci teknik, gizli ısı depolamadır. Faz değişim malzemesi (PCM) olarak da adlandırılan bir termal enerji depolama malzemesi ile depolama sağlanır. Faz değişim sıcaklığı olarak bilinen belirli bir sıcaklıkta fazı değiştirir. Başka bir deyişle, PCM'deki kimyasal bağlar, sıcaklık faz değiştirme sıcaklığının üstüne çıktığında kopmaya başlayacaktır. O sırada malzeme erimeye başlar. Sıcaklık düştükçe malzeme ekzotermik reaksiyon nedeniyle katılaşır. Üçüncü termal depolama tekniği, termal enerjinin kimyasal reaksiyonlarla depolanmasına dayanan termo-kimyasal enerji depolamadır. Termo-kimyasal enerji depolamanın ana mekanizması, iki veya daha fazla bileşenin bir kimyasal bileşikte birleştirildiğinde, bu bileşiğin daha sonra ısı ile bölünmesidir. Bölünmüş bileşenler daha sonra ayrı olarak saklanır. Talep durumunda, bileşenler kimyasal bir bileşik halinde birleştirilir ve ısı açığa çıkar (Gustavsson, 2016).

## 3. MİKRO ŞEBEKE TASARIMINDA VE İŞLETİLMESİNDE KARŞILAŞILAN SORUNLAR VE ÇÖZÜM ÖNERİLERİ

Bu çalışmada, enerji üretimindeki 3 önemli operasyonel değişken ele alınmıştır: fazladan üretilen elektrik, $CO_2$ emisyonu ve su temini kıtlığı. Bu bölümde bu sorunlar ve çözüm stratejileri ayrıntılı olarak açıklanmıştır.

### 3.1. Fazladan Üretilen Elektrik Sorunu ve Olası Çözümleri

Fazladan elektrik üretimi önemli bir ekonomik sorundur. Bunun nedeni, talepten fazla üretilen elektriğin depolanması ile ilgili zorluklardır. Piller bu soruna bir çözüm olabilir. Bununla birlikte, enerji sistemlerine pil eklemek maliyeti arttırır ve bu da karı azaltabilir. Alternatif olarak, fazla elektriği ana şebekeye veya enerji sistemine yakın bir yere zamanlama ve talebin uyuşması durumunda satmak mantıklı bir çözüm olabilir.

Elektroliz, fazla elektriği kullanmak için makul bir çözüm olabilir. Elektroliz, büyük miktarlarda elektrik tüketerek suyu hidrojen ve oksijene dönüştürür. Hidrojen, birçok endüstriyel uygulama için en çok tercih edilen enerji kaynağıdır. Hidrojen çeşitli oranlarda CO ile karıştırılarak sentez gazı üretiminde kullanılabilir. Bu da dizel ve benzin üretimi için bir hammaddedir. Üretilen hidrojen depolanmayabilir. Bunun yerine, depolama maliyetini düşürmek için üretildiği gibi tüketilmesi tercih edilir. Amonyak sentezi hidrojeni tüketmek için iyi bir alternatif olabilir. Haber-Bosch prosesi, hidrojen ve nitrojenden amonyak üreten termokatalitik bir prosestir (Singh ve diğ., 2017).

Birkaç elektroliz türü vardır. İlki, hücreleri bir tür sıvı elektrolit olan alkalin elektrolizdir. İkincisi, PEM (proton değişim membranı) elektrolizidir. Alkaline göre daha yüksek akımla çalışır.

Ayrıca, iyon transfer ortamı olarak bir polimer elektrolit asit membranı kullanır. Diğer bir teknoloji, katı oksit elektroliz hücresi (SOEC) elektrolizidir. Bahsedilen diğer tipler için uygun olan sıcaklıklardan çok daha yüksek sıcaklıklarda çalışır. Enerji sistemindeki hem karbondioksiti hem de fazla elektriği azaltmak için, yüksek sıcaklıkta SOEC elektroliz sistemi ve katalitik metanasyon sistemleri birleştirilebilir. Oluşacak kombine sistem, suyu hidrojene dönüştürüp ardından karbondioksitle birleştirerek sentez gazı oluşturabilir. Öte yandan, SOEC sistemi Fischer-Tropsch prosesi ile birleştirilebilir. Bu sistem, enerji sistemindeki fazla elektriği tüketerek hidrokarbon sentezler (Giglio ve diğ., 2015; Kaczur ve diğ., 2018).

### 3.2. $CO_2$ Emisyonu Sorunu ve Olası Çözümleri

2018'de yayınlanan Hükümetler Arası İklim Değişikliği Paneli (IPCC) raporuna göre, sera gazı salınımları küresel ısınmayı en fazla 1,5 °C arttıracak şekilde sınırlandırılmalıdır. Ayrıca, karbondioksit, su buharı, azot oksit ve metan gibi küresel sera gazı emisyonlarının net değerinin, 2050 yılına kadar sıfıra indirilmesi gerekmektedir. Sera gazı emisyonlarını azaltmak oldukça zor bir hedeftir. Örneğin, tarımsal salınımların sıfıra indirilmesi güncel teknoloji ile çok zordur. Bunun yerine; elektrik üretimi, endüstri ve ulaşım sistemleri gibi çeşitli alanlardan kaynaklanan salınımların azaltılması daha önceliklidir.

Üretilen karbondioksit temelde iki şekilde uzaklaştırılabilir: son kullanım için karbon tutma (CCU) ve depolama için karbon tutma (CCS). CCS'de endüstriyel işletmelerin baca gazlarından çıkan karbondioksit toplanır ve daha sonra yer altı rezervuarlarında depolanır. CCU'da ise karbondioksit işletmeden uzaklaştırıldıktan sonra toplanarak başka bir işlemde kullanılmaktadır. Örneğin, bir işlemden çıkarılan karbondioksit, ticari olarak değerli bir kimyasal olan formik aside elektrokimyasal olarak indirgenebilir (Aldaco ve diğ., 2019).

Karbondioksit, madencilik endüstrisinde kayaların kırılması için de kullanılmaktadır. Örneğin, killi şeyl (şist) oluşumları dünya genelinde yaygındır. Bu tür şist oluşumlarını kırmak için su bazlı olmayan sıvılar gereklidir; çünkü su bazlı kırılma, şişmeyi arttırır ve şistin gaz depolama kapasitesini azaltır. $CO_2$, küçük şişme etkisi ve düşük kırılma basınçlarında bile uzun ve dar çatlaklar oluşturma yeteneği ve çok hızlı temizleme özelliği sebebiyle sıklıkla tercih edilmektedir (C. Zhang ve diğ., 2017).

Metan ve metanol sentezi, hammadde olarak karbondioksit ve hidrojene ihtiyaç duyulan diğer uygulamalardır. Metan, çok aşamalı bir prosesle üretilir. Bu proses, Sabatier prosesi adı verilen bir tür metanasyon-elektroliz kombine işlemiyle başlar. Sabatier işleminde, hidrojenin karbondioksit ile reaksiyona girmesine izin verilerek metan ve su elde edilir (Gustavsson, 2016).

Metanol sentezi, kimyasal uygulamalar arasında hidrojeni tüketen ikinci yaygın prosestir. Metanol, buhar reformasyonu ve kömür gazlaştırma gibi sentez gazı proseslerine dayalı olarak üretilir. Hafif doğal gaz gibi hafif hammaddeden kükürt uzaklaştırılmaktadır, daha sonra buhar ve oksijen ile reformasyon yapılmaktadır. Elde edilen sentetik gaz sıkıştırılmaktadır ve metanol-sentez ünitesine aktarılmaktadır. Sentez gazını metanole dönüştüren reaksiyonlar oldukça ekzotermiktir. Bu nedenle; metanol sentezinin verimliliği, ısı geri kazanımı ve yönetiminin verimliliği ile doğrudan ilişkilidir (Giglio ve diğ., 2015).

### 3.3. Su Kaynaklarının Yetersizliği Sorunu ve Olası Çözümleri

Fazla elektriği azaltmak için, bu çalışmanın önceki bölümünde suyun elektrolizi önerilmiştir. Elektroliz işlemi alternatif olarak düşünülürken, mikro şebekedeki su kaynağının yeterli olup olmadığı dikkate alınmalıdır. Küresel perspektiften bakıldığında, su kaynakları her geçen gün azalmaktadır. Su kıtlığı, kısa sürede tüm dünyada kıtlık ve kuraklık getirecek bir çevre sorunu haline gelmiştir. Su arzının azalması, su fiyatını etkileyebilir ve ülkeler arasında bir rekabet başlatabilir. Bu büyük sorunun üstesinden gelmek için, içme suyu üreten teknolojiler kullanılmalıdır veya mevcut teknolojiler kurulacak büyük ölçekli tesisler için iyileştirilmelidir.

Suyu tuzdan arındırma, günümüz adalarında yaygın olarak kullanılan bir tür su arıtma tekniğidir. Bununla birlikte, su temini kıtlığı tehdidi düşünüldüğünde, bu teknik dünyanın her yerine yayılabilir. Bu nedenle, suyun tuzdan arındırma teknolojisi son yıllarda çok daha fazla ilgi görmektedir. Suyu tuzdan arındırma teknolojisi, tatlı su elde etmek için tuzlu sudan çözünmüş tuzları ve mineralleri ayırma işlemi olarak tanımlanır. Çeşitli tuzdan arındırma teknolojileri arasında iki temel teknik popülerdir: buharlaşmaya dayalı işlemler (faz değişimi) ve membranlara dayalı işlemler (faz dışı değişim). Buhar üretmek için deniz suyu veya acı su gibi besleme suyunun işletme basıncında kaynama noktasına kadar ısıtılmasıyla buharlaşma bazlı tuzdan arındırma yapılır. Daha sonra bu buhar, tatlı su elde etmek için bir ünitede yoğunlaştırılır. Buharlaşmaya dayalı tuzdan arındırma, çok aşamalıdır. Damıtma, termal buhar sıkıştırması veya mekanik buhar sıkıştırma gibi prosesleri içermektedir. Membran bazlı tuzdan arındırma tekniği üç farklı yolla yapılabilir. Bunlar; elektrik enerjisi kullanan ters ozmos prosesi, nanofiltrasyon ve iyonik membranların kullanıldığı ters elektrodiyaliz prosesidir.

Ters ozmos yönteminde su, yarı geçirgen ince membranlardan yüksek basınçta geçirilerek yüksek tuz konsantrasyonu düşürülür. Dünya genelinde kurulu su arıtma tesislerinin % 50'den fazlası, düşük enerji maliyeti ve elektrik tüketimi sebebiyle ters ozmos yöntemini kullanmaktadır. Ters ozmos işlemleri, modüler tasarımları nedeniyle genellikle küçük ve orta kapasiteli tuzdan arındırma tesislerinde gerçekleştirilir. Ters ozmos prosesleri, suyun tuzdan arındırılması için gerekli elektriği üretmek için özellikle fotovoltaik hücreler olmak üzere yenilenebilir jeneratörlerle ile birleştirilmektedir (Calise ve diğ., 2019). Şebekelerde üretilen fazla elektriğin, bu tip su arındırma süreçlerine dahil edilmesi de etkili bir yöntem olarak düşünülebilir.

## 4. KARIŞIK (KARMA) TAM SAYILI DOĞRUSAL OLMAYAN PROGRAMLAMA

Bir işletmenin tasarımı, konumunun seçimi ve üretiminin planlaması gibi konuları içeren çok sayıda optimizasyon problemi, sadece sürekli değil, aynı zamanda kesikli (ayrık) olan karar değişkenleri ile modellenir. Örneğin, yeni bir ekipman parçasının kurulup kurulmadığına dair bir kararı açıklayan değişkenler, "0-1" değişkeni veya 'ikili değişken' tanımlanır.

Bu tür optimizasyon problemleri, süreksiz (kesikli) değişkenlerinin, amaç fonksiyonlarının ve kısıtların matematiksel yapısına göre sınıflandırılabilir. Karma tam sayı programlama (MIP), amaç fonksiyonunun sürekli ve kesikli (tam sayı) olan iki değişken setine bağlı olduğu bir programlama olarak tanımlanır. Bir optimizasyon probleminin sadece tam sayı değişkenleri varsa, buna tam sayı programlama (IP) problemi denir. Ayrıca, IP'deki bu tam sayı değişkenleri ikili değişkenler ise, yani karar değişkenleri olarak "0-1" içeriyorsa bu problemlere, ikili tam sayı programlama (BIP) denir.

MIP problemleri, doğrusal veya doğrusal olmayan amaç fonksiyonlarına ve kısıtlarına göre iki ana sınıfa ayrılır. MIP'nin yalnızca doğrusal amaç fonksiyonları ve kısıtları varsa, buna karma tam sayı doğrusal programlama (MILP) problemleri denir. Öte yandan, gerçek operasyonel koşullarda optimizasyon problemleri her zaman doğrusal olmayabilir. Bunun yerine, bazı doğrusal olmayan amaç fonksiyonları veya kısıtlar içerirler. Bu durumda, problem karma tam sayılı doğrusal olmayan programlama (MINLP) olarak adlandırılır. MINLP genellikle hesaplama açısından zor bir optimizasyon problemidir; ama son yirmi yıldır hesaplama alanında çok önemli gelişmeler meydana gelmiştir. Bu nedenle, günümüzde MINLP problemleri için çeşitli etkili çözücüler mevcuttur.

Dışbükey karma tam sayılı doğrusal olmayan programlamada en sık tercih edilen algoritmalar dal-sınır, doğrusal olmayan programlama tabanlı dal-sınır, genişletilmiş kesme düzlemi, dışsal yaklaşıklama, genellendirilmiş Benders ayrıştırması ve genişletilmiş destek hiper düzlemidir.

Basit bir MINLP problemi şu şekilde formüle edilebilir:

$$\underset{x,y}{Min} c_1^T . x + c_2^T . y \qquad (1)$$

$$N = g_i(x, y) \leq 0 \qquad (2)$$

$$L = Ax + By \leq b \qquad (3)$$

$$c \leq y_i \leq d, \forall i = 1, 2, \ldots m \qquad (4)$$

Denklem 1-4'te, doğrusal olmayan (*N*), doğrusal (*L*) ve ikili değişken (*y*) kısıtları verilmiştir. $y_i$ ayrık değişken, *x* sürekli değişken, $c_1$ ve $c_2$ matris katsayıları, *b*, *c* and *d* gerçel sayılardır. Dal-sınır algoritması, orijinal problemin tam sayı kısıtlarını gevşeterek ve ardından sürekli konveks doğrusal olmayan programlama gevşetmelerini çözerek karma tam sayılı doğrusal olmayan programlama problemini çözer. Problemin sürekli gevşetmesini çözmek, geçerli bir alt sınır oluşturan bir çözüm ($x^k$, $y^k$) verir. $y^k$'nın tüm bileşenleri tam sayı değerleri alırsa, karma tam sayılı doğrusal olmayan programlama probleminin optimum çözümü var demektir. Değilse, sürekli gevşetme, gevşetilmiş probleme $y_i \leq y_i^k$ ve $y_i \geq y_i^k$ kısıtları eklenerek problem iki yeni doğrusal olmayan alt probleme ayrılır. Yeni alt problemler çözülerek yeni bir alt sınır oluşturulabilir. Alt problemlerden biri bir tam sayı çözümü üretirse, bu aynı zamanda geçerli bir üst sınırdır. Optimum çözüm için arama prosedürü, düğümlerin alt problemleri göstererek kendi ana düğümlerine bağlandığı bir ağaçla gösterilmektedir. Düğümlerden biri tam sayı çözümü vermezse, iki yeni alt problem oluşturan iki yeni düğüme bölünür. Düğümlerden biri üst sınırdan daha kötü bir optimuma sahipse veya alt problem uygulanabilir değilse, o zaman düğüm kesilebilir çünkü orada optimum bir çözüm bulunamaz. Dallanma ve sınır algoritması her düğümde dışbükey doğrusal olmayan programlama problemlerinde kullanılıyorsa, buna doğrusal olmayan programlama tabanlı dal-sınır algoritması denir. MINLP'de, her gevşetilmiş alt problem dışbükey değilse, birçok NLP algoritması küresel bir çözüme yakınsamaya ulaşabilir. Dal-kesim ve dal-fiyat algoritmaları dal-sınır algoritmasının farklı varyasyonlarıdır (Quesada ve Grossmann, 1992).

Genişletilmiş kesme düzlemi (ECP) algoritması, doğrusal olmayan kısıtları doğrusal olarak iyileştiren iterasyonlarla çok yüzlü dışsal yaklaşıklama geliştiren bir algoritmadır. Genişletilmiş destek hiper düzlem (ESH) algoritması, dışbükey tam sayılı doğrusal olmayan programlama problemlerini çözmek için kullanılan diğer bir algoritmadır. ESH algoritması, deneme çözümleri bulabilmek için ECP algoritmasıyla aynı prensipte çalışır. Bununla birlikte, çok yüzlü dışsal yaklaşıklama oluşturmak için farklı metotlar kullanır. ECP algoritmasında çok yüzlü dışsal yaklaşıklama oluşturmak için şekillendirilmiş kesme düzlemleri yeterince sıkı olmayabilir. Bununla birlikte, ESH algoritması, çok yüzlü bir dış yaklaşım oluşturmak için her yinelemede destekleyici hiper düzlemler elde edebilir (Westerlund ve Pettersson, 1995).

Dışsal yaklaşıklama (OA), bir dizi karma tam sayılı doğrusal programlama ve doğrusal olmayan programlama alt problemini çözen orijinal problemin optimum çözümünü veren bir ayrıştırma yöntemidir. Hem ECP hem de ESH'ye benzer şekilde çalışır. İlaveten, dışsal yaklaşıklama doğrusal olmayan uygulanabilir bölgenin çok yüzlü dış yaklaşımını yinelemeli olarak iyileştirir. Dışsal yaklaşıklama, $y^k$ tam sayı kombinasyonunu seçmek için yalnızca çok yüzlü yaklaşımını kullanırken, $x^k$ değerini seçerken bir dışbükey doğrusal olmayan alt problemi çözer (Duran ve Grossmann, 1986).

Genellendirilmiş Benders ayrıştırması (GBD) algoritması, karma tam sayılı doğrusal programlama problemlerini çözmek için geliştirilen bir bölümleme prosedürüdür. Bu algoritma, dışsal yaklaşıklama algoritması ile yakından ilişkilidir. Aralarındaki temel fark, ana problemin türetilmesidir. Genellendirilmiş Benders ayrıştırmasında, ana problem tam sayı değişkenleri tarafından tanımlanan uzaya yansıtıldığından, yalnızca tam sayı değişkenleri ifade edilir (Geoffrion, 1972).

AlphaECP, Antigone, AOA, BARON, BONMIN, Couenne, DICOPT, Juniper, KNITRO, LINDO, Minotaur, Muriqui, Pavito, SBB, SCIP ve SHOT karma tam sayılı doğrusal olmayan programlama çözücülerinden bazılarıdır. Çözücülerin çoğunluğu problemleri çözerken tek bir algoritma kullanmazlar. Bunun yerine çeşitli teknikleri birleştirerek çözümü iyileştirmeyi hedeflerler. Optimizasyon çözücülerinin çoğu, optimizasyon için AIMMS, AMPL ve GAMS gibi iyi organize edilmiş modelleme ortamına bağlıdır. Ayrıca Python ve Julia'da da optimizasyon modellemesi yapılmasına yükselen bir talep vardır. JuMP, Julia'da bulunan bir optimizasyon modelleme ortamı iken; Pyomo, Python'daki modelleme ortamıdır. Bazı MINLP çözücüleri, MATLAB'a ara yüz de sağlar. OPTI Toolbox, MATLAB'te çeşitli çözücülere ulaşmak için kullanılır. Bu çalışmada, BARON, SCIP ve DICOPT çözücülerinin özelliklerine yer verilmiştir.

BARON, dal-azaltma (branch and reduce) algoritması kullanarak global çözüm verir. Dal-sınır düğümlerindeki doğrusal programlama gevşemelerini çözer. Ayrıca BARON, karma tam sayılı doğrusal programlama gevşemelerini ve doğrusal olmayan gevşemeleri çözer. Dışbükey olmayan problemleri, bir uzamsal dal-sınır tekniği ile dışbükey eksik tahmin ediciler ve içbükey aşırı tahmin edicileri birleştirerek çözer. Çözücü, dışbükey olmayan işlevleri dışbükey veya içbükey gevşemelere sahip daha basit parçalara ayırmak için otomatik yeniden formülasyonlar ve dışbükeylik tanımlaması yapar. Yeniden formülasyonlar, daha sıkı çok yüzlü dış yaklaşımlar sağlayabilir. BARON ayrıca arama alanını düşürmek için gelişmiş sınır sıkma ve menzil azaltma yöntemlerini kullanır. Buna ek olarak, yerel arama tekniklerini birincil sezgisel tarama olarak gerçekleştirir. BARON, doğrusal olmayan programlama alt problemlerini çözmek için iç nokta iyileştirici (IPOPT), FilterSD veya FilterSQP'ı kullanır. BARON ayrıca GAMS'te mevcut herhangi bir doğrusal olmayan çözücüyü çalıştırabilir.

SCIP (Çözme Kısıtlı Tam sayı Programları) çözücüsü Berlin Zuse Enstitüsü'nde geliştirilmiştir. Akademik amaçlı kullanıcılar için ücretsizdir. Dal-kes-fiyat algoritmasını kullanarak kısıt tam sayılı ve karma tam sayılı programlama problemleri çözer. SCIP, çok yüzlü dış yaklaşımlar ve bir uzaysal dal ve sınır algoritması kullanarak dışbükey ve dışbükey olmayan MINLP problemini çözmek için geliştirilmiştir. Çözücü, güçlü dual sınırlar oluşturmak için doğrusal programlama gevşemelerini ve kesme düzlemlerini kullanır. Çözücüde ayrıca çeşitli ilkel buluşsal yöntemler ve sınır sıkma yöntemleri de kullanılır. SCIP, doğrusal programlama alt problemlerini çözmek için SoPlex'i kullanır. Bununla birlikte, varsa CLP, CPLEX, Gurobi, Mosek veya XPress'i de kullanabilir. Ayrıca çözücü, doğrusal olmayan programlama alt problemlerini çözmek için IPOPT kullanır (Anand ve diğ., 2017; Kronqvist ve diğ., 2019; Melo ve diğ., 2020).

DICOPT, doğrusal ikili veya tamsayı değişkenleri ve doğrusal ve doğrusal olmayan sürekli değişkenleri içeren karma tamsayılı doğrusal olmayan programlama problemlerini çözer. MINLP optimizasyon problemlerinin modellemesi ve çözümü henüz doğrusal, tamsayılı veya doğrusal olmayan programlama modellemesi ile elde edilen olgunluk ve güvenilirlik aşamasına ulaşmamıştır. DICOPT, Carnegie Mellon Üniversitesi Mühendislik Tasarım Araştırma Merkezi'nde J. Viswanathan ve Ignacio E. Grossmann tarafından geliştirilmiştir. Program, dışsal yaklaşıklama algoritmasının uzantılarına dayanmaktadır. DICOPT içindeki MINLP algoritması, bir dizi NLP ve MIP alt problemini çözer. Bu alt problemler, GAMS altında çalışan herhangi bir NLP (Doğrusal Olmayan Programlama) veya MIP (Karma Tamsayılı Programlama) çözücü kullanılarak çözülebilir. Performans, büyük ölçüde seçilen alt çözücülerin seçimine bağlı olacaktır. Algoritma, dışbükey olmayanları işlemek için gerekli niteliklere sahip olsa da, küresel optimumu mutlaka elde etmek zorunda değildir (GAMS, 2021).

## 5. OPTİMİZASYONUN MODELİNİN OLUŞTURULMASI

Enerji yönetim sistemi tasarım problemini çözmek için bir optimizasyon modelinin formüle edilmesi gerekmektedir. Optimizasyon modelleri, optimize edilmesi gereken değeri hesaplayan bir amaç fonksiyonunu, bileşen modellerini, çalışma koşullarını ve diğer gereksinimleri ifade eden eşitlik ve eşitsizlik kısıtlarını içerir.

### 5.1. Amaç Fonksiyonu ve Kısıtlar

Amaç fonksiyonları farklı şekillerde yazılabilir. En bilinen amaç fonksiyonu, ekonomik amaç fonksiyonudur. Bu fonksiyon, kârı veya net bugünkü değeri maksimize etmek ve işletme maliyetini veya toplam maliyeti en aza indirmek vb. şeklinde yazılabilir. Boyuta bağlı olan maliyet fonksiyonu, denklem 5'teki gibi doğrusal olabilir veya denklem 6'daki gibi ikinci dereceden doğrusal olmayan bir fonksiyon olabilir. Denklemlerdeki $P_r$ ekipmanın anma gücünü ve $a, b, c, d, e, f$ denklemlere özgü katsayıları göstermektedir.

$$\textit{Kurulum Maliyeti} = a + b.P_r \tag{5}$$

$$\textit{Kurulum Maliyeti} = d + e.P_r + f.P_r^2 \tag{6}$$

Doğrusal maliyet fonksiyonu her zaman artan boyutla doğrusal olarak artar, oysa doğrusal olmayan maliyet işlevi boyut arttıkça azalan ivme ile artar. Amaç fonksiyonu, enerji sistemindeki yakıt tüketimi ile ilgili karbondioksit salınımlarının en aza indirilmesi gibi çevresel kısıtlara dayalı olarak da yazılabilir.

Kurulan modelin amacı, mikro şebekeye kurulacak optimum ekipman sayısını bulmak ise, modele aday ekipmanlar tanıtılmalıdır. Ekipman seçimleri yapılırken, ikili (binary) değişken ('0-1') ekipmanın kurulup kurulmayacağı kararını verir. Eğer optimizasyon sonunda modelin seçtiği karar değişkeni '0' ise o ekipmanın kurulmaması optimum tasarıma uygundur. Öte yandan, eğer '1' ise o ekipmanın kurulması optimum tasarıma uygundur. Benzer şekilde; üretim planlaması üzerine bir optimizasyon yapılıyorsa, ikili değişkenin '0' olduğu zaman diliminde kurulu ekipmanın çalıştırılmaması, '1' ise çalıştırılması optimum planlamaya uygundur. Ekipmanlara ait model denklemleri bir sonraki bölümde açıklanacaktır.

Ekipmanların model denklemleri oluşturulduktan sonra, sistemin enerji denkliği yazılmalıdır. Müşteri enerji talebini de içeren enerji denkliği, optimizasyon modelinde eşitlik kısıtı olarak tanımlanmıştır. Enerji denkliği oluşturulurken, enerji sisteminin her bir bileşeni enerji girdisi ya da çıktısı şeklinde gösterilir. Bir enerji üretim şebekesine ait enerji akış diyagramı Şekil 1'de gösterilmiştir.

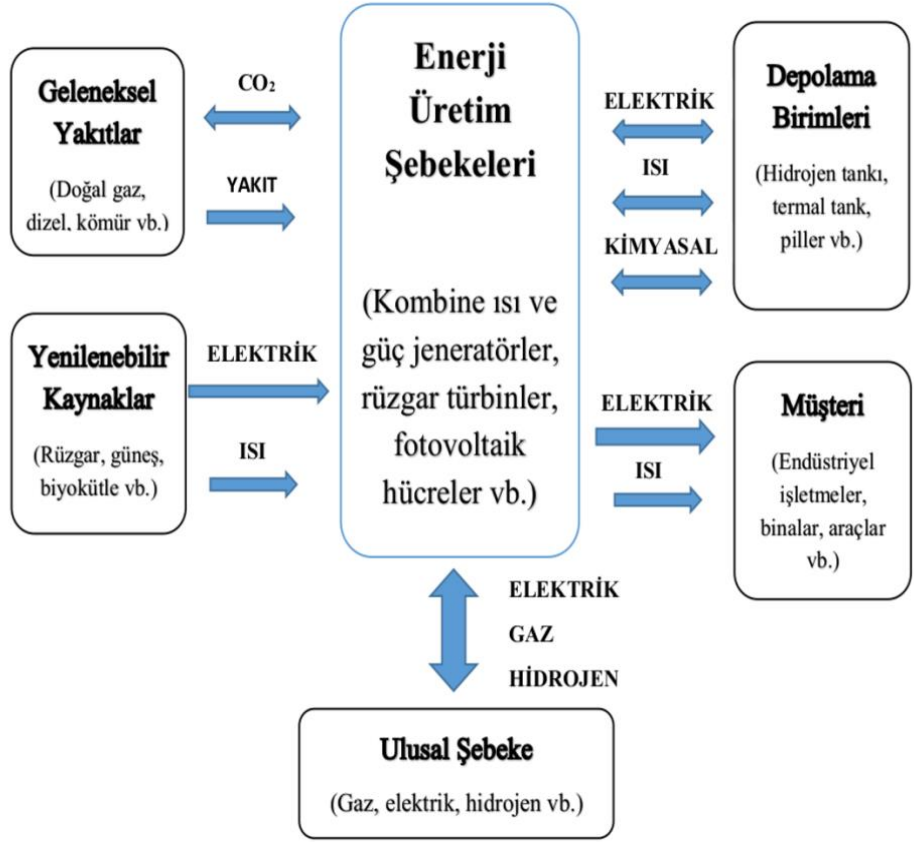

**Şekil 1.** Enerji üretim şebekesinde basit bir enerji akış şeması.

Mikro şebekenin çalışma koşulları, eşitsizlik veya eşitlik kısıtları olarak modele eklenmelidir. Bu kısıtlardan bazıları sıralanmıştır:

1) Maksimum seçilebilecek ekipman sayısı: Mikro şebekenin kurulacağı alan, toplam ekipman sayısını sınırlandırır.

2) Kapasite sınırları: Her bir ekipmanın anma gücü, üreticisinin sağladığı minimum ve maksimum güç aralığında olmalıdır.

3) Çalışma sınırları: Her bir ekipmanın çalışma yükü, ekipmanın nominal yükünü aşamaz.

4) Bir ekipmanın çalışmasının hızlı dalgalanmasını önlemek için iki ardışık zaman anındaki yük farkı sınırlanmalıdır.

5) $CO_2$ salınımları: Salınım limitleri, hükümet tarafından belirlenen değerler ile sınırlandırılmalıdır.

6) Kombine ısı ve güç jeneratörleri için sıcaklık, jeneratör üreticisinin verdiği maksimum değeri aşmamalıdır.

7) Depolama sınırları: Bir depolama cihazındaki enerji miktarı, depolama kapasitesini aşmamalıdır.

8) Enerjinin depolanmasının sürekliliği: Depolama cihazının herhangi bir gündeki son durumu, sonraki günkü başlangıç durumu ile aynı olmalıdır (Abo-Elyousr ve Elnozahy, 2018; Aghaei ve Alizadeh, 2013; Alvarado-Barrios ve diğ., 2020).

### 5.2. Mikro Şebeke Bileşenlerinin Model Denklemleri

Bu bölümde, jeneratörlerin model denklemleri ve depolama cihazlarının şarj/deşarj modeli tanıtılmaktadır.

#### 5.2.1. Kombine Isı ve Güç (CHP) Jeneratörleri

CHP sistemleri, enerji verimliliğini artırmak için aynı anda güç ve ısı üretebilir. CHP üniteleri geri basınç ve ekstraksiyon-yoğunlaştırma üniteleri olarak sınıflandırılır. Geri basınç ünitelerinde, üretilen ısı güç çıktısıyla doğrusal bir ilişki içindedir. Isı/güç oranı sabittir, bu durum sınırlı bir operasyon esnekliği sağlar. Ekstraksiyon-yoğunlaştırma üniteleri için operasyon bölgesi düzensiz bir dörtgen olarak tanımlanmaktadır. Isı/güç oranı değişkendir. Bu nedenle, ekstraksiyon-yoğunlaştırma ünitelerinin uygulanabilir operasyon bölgesi (FOR), geri basınç ünitelerininkinden çok daha büyüktür. Güç ve ısı arasındaki ilişkiyi açıklayan bu bölgeye ait kısıtlar modele eklenmelidir (Guo ve diğ., 2019; Ko ve Kim, 2019). CHP sistemlerinin işletme maliyeti, yakıt tüketim maliyetine ve bakım-onarım maliyetlerine bağlıdır. Yakıt tüketim maliyeti, üretilen ısı ve güce bağlı doğrusal olmayan bir fonksiyondur (Roy ve diğ., 2014):

$$M(H_{chp}, P_{chp}) = a_{chp} + b_{chp}P_{chp} + c_{chp}P_{chp}^2 + d_{chp}H_{chp} + e_{chp}H_{chp}^2 + f_{chp}H_{chp}P_{chp} \tag{7}$$

$M(H_{chp}, P_{chp})$, saatlik yakıt tüketim maliyetidir. $H_{chp}$, çalıştırıldığı saat diliminde chp sisteminin ürettiği ısı; $P_{chp}$ ise ürettiği güçtür. $a_{chp}, b_{chp}, c_{chp}, d_{chp}, e_{chp}$ ve $f_{chp}$ yakıt tüketim fonksiyonun sabit terimli katsayılarıdır.

CHP sistemlerinin saatlik toplam gaz emisyonu, üretilen güce bağlı doğrusal bir fonksiyonla hesaplanır. Denklem 8'de, $CO_2$ emisyonunun toplam gaz salınımının %99'u olduğu kabul edilmiştir (Nwulu, 2020). Karbondioksit emisyonu fonksiyonu:

$$CO2_{chp} \ (ton/saat) = 0,99(\emptyset_{chp} + \mu_{chp})P_{chp} \tag{8}$$

$\emptyset_{chp}$ ve $\mu_{chp}$ fonksiyonunun sabit terimli katsayılarıdır.

#### 5.2.2. Konvansiyonel Enerji Jeneratörleri

Konvansiyonel enerji jeneratörleri; kömür, doğal gaz ve petrol gibi yakıtları kullanarak elektrik üretir. Kurulum ve işletme maliyetleri ne tip yakıt kullanıldığına bağlıdır. Yakıt tipi karbondioksit salınımlarını da etkiler. Konvansiyonel jeneratörler için yakıt tüketim maliyeti denklem 9'da verilmiştir:

$$M_{konv} = c_{konv}(P_{konv}^2) + b_{konv}(P_{konv}) + a_{konv} + \left| d_{konv} \sin\left(e_{konv}(P_{konv}^{min} - P_{konv})\right) \right| \tag{9}$$

$M_{konv}$, saatlik yakıt tüketim maliyetidir. $P_{konv}$, jeneratörün çalıştırıldığı saat dilimindeki ürettiği güç ve $P_{konv}^{min}$ ise jeneratörün ürettiği minimum güçtür. $a_{konv}, b_{konv}, c_{konv}, d_{konv}$ ve $e_{konv}$ fonksiyonun sabit terimli katsayılarıdır (Mohammadi-Ivatloo ve diğ., 2013).

Konvansiyonel güç jeneratörlerinin karbondioksit emisyonu, kullandıkları yakıt tipiyle doğrudan orantılıdır. Karbondioksit emisyonu fonksiyonu üretilen güce bağlı ikinci dereceden polinom fonksiyondur:

$$CO2_{konv}\left(\frac{ton}{saat}\right) = ef_{konv}(h_{konv}P_{konv}^2 + g_{konv}P_{konv} + f_{konv}) \quad (10)$$

Denklem 10'da $ef_{konv}$, yakıt tüketimine göre değişen sabit bir katsayıdır. $h_{konv}$, $g_{konv}$ ve $f_{konv}$ polinom fonksiyonun sabit terimli katsayılarıdır (Y. Zhang ve diğ., 2013).

### 5.2.3. Fotovoltaik Hücreler

Fotovoltaik hücrelerin güç çıktısı aşağıdaki şekilde hesaplanabilir:

$$P_{FV}(t) = P_{FVref} \cdot \frac{S(t)}{S_{ref}}[1 + B_{ref}(T_c(t) - T_{ref})] \quad (11)$$

$$T_c(t) = T_a(t) + \frac{(NOCT - 20)}{800} \cdot S \quad (12)$$

Denklem 11-12'de *t* anlık zamandır, $P_{FV}(t)$ fotovoltaik hücrenin o anda ürettiği güç çıktısı, $P_{FVref}$ referans alınan hücrenin güç çıktısı, *S(t)* anlık zamandaki radyasyondur, $T_c$ fotovoltaik hücrenin yüzey sıcaklığıdır, $T_a$ ortam sıcaklığıdır (°C), *NOCT* nominal işletim hücre sıcaklığıdır. *NOCT*, fotovoltaik modülün 1000 W/m² güneş ışınımı ve 20 °C ortam sıcaklığı altında çalıştığı hücre sıcaklığıdır. Genellikle 42º C ile 46 ºC arasındadır. $B_{ref}$, silikon hücreler için °C başına 0,004 ila 0,006 arasında değişen sıcaklık verimlilik katsayısıdır. *S*, mW/cm² cinsinden güneş radyasyonu seviyesidir (Abo-Elyousr ve Elnozahy, 2018; Honsberg ve Bowden, 2021).

Mikro şebeke optimizasyonunda, fotovoltaik hücrede üretilen güç çıktısı, anma gücüne ve o zaman dilimindeki güneş radyasyonuna bağlı oranlarla ifade edilebilir:

$$P_{FV}(R) = \begin{cases} P_{nominal}\left(\dfrac{R^2}{R_C R_{STD}}\right), & 0 \leq R \leq R_C \\ P_{nominal}\left(\dfrac{R}{R_{STD}}\right), & R_C < R < R_{STD} \\ P_{nominal}, & R_{STD} \leq R \end{cases} \quad (13)$$

Denklem 13'te $P_{nominal}$, fotovoltaik hücrenin anma gücüdür. *R*, radyasyondur (W/m²). $R_{STD}$, standart koşullardaki güneş radyasyonudur ve genellikle 1000 W/m² olarak alınır. $R_C$, belirli bir radyasyon noktasını temsil eder ve genellikle 150 W/m² olarak hesaplamalara dahil edilir (Ben Hmida ve diğ., 2019; Morshed ve diğ., 2018).

### 5.2.2. Rüzgar Türbinleri

Rüzgar türbinlerinin güç çıktısı aşağıdaki şekilde ifade edilir:

$$P_{RT} = \frac{1}{2}\rho A C_p V^3 \quad (14)$$

Denklemde $\rho$ havanın yoğunluğu, *A* süpürülen alan, $C_p$ güç katsayısı ve *V* rüzgar hızıdır. (Abo-Elyousr ve Elnozahy, 2018; Amano, 2017).

Mikro şebeke optimizasyonunda, rüzgar türbininde üretilen güç çıktısı, anma gücüne ve o zaman dilimindeki rüzgar hızına bağlı oranlarla ifade edilebilir:

$$P_{RT} = \begin{cases} 0 & , v < v_{cut,in} \\ 0 & , v > v_{cut,out} \\ P_{nominal} & , v_{nominal} \leq v \leq v_{cut,out} \\ P_{nominal}\left(\dfrac{v - v_{cut,in}}{v_{nominal} - v_{cut,in}}\right) & , v_{cut,in} \leq v \leq v_{nominal} \end{cases} \quad (15)$$

Denklem 15'te, $v_{nominal}$ nominal rüzgar hızıdır. Cut-in hızının altında ve Cut-out hızının üstündeki rüzgar hızlarında üretilen güç çıktısı sıfırdır (Dolatabadi ve diğ., 2017).

### 5.2.3. Depolama Birimleri

Bataryalar, hidrojen tankları ve termal tanklar gibi depolama cihazları, şarj ve deşarj modlarında çalışır. Depolama cihazlarının şarj/deşarj modeli şu şekilde yazılabilir:

$$X_t = X_{(t-1)} \cdot \eta + \Delta t \cdot (\eta_{ch} X_t^{ch} - X_t^{dch}/\eta_{dch}) \quad (16)$$

Denklemde $\eta$ bekleme verimliliği, $\eta_{ch}$ şarj verimliliği, $\eta_{dch}$ deşarj verimliliği, $X_t$ depolama cihazında $t$ zamanındaki enerjidir, $X_{ch}$ ve $X_{dch}$ şarj/deşarj hızıdır (Rancilio ve diğ., 2019).

## 6. SONUÇ

Bu çalışmada, fosil ve yenilenebilir enerji kaynakları içeren mikro şebekelerin optimum tasarımı ve üretim planlaması incelenmiştir. Fazladan elektrik üretimi, $CO_2$ emisyonu, su kıtlığı gibi problemlere yönelik olası çözümler ele alınmıştır. Bu problemler için önerilen optimizasyon metodu, literatürde henüz çok da fazla kullanılmamış olan karma tam sayılı doğrusal olmayan programlamadır (MINLP). Son olarak, optimizasyon modeli oluşturulurken dikkate alınması gereken kısıtlar ve her bir ekipmana ait model denklemleri verilmiştir.

## 7. TEŞEKKÜR

## KAYNAKLAR